\documentclass[12pt,a4paper]{article}
\title{Invariant distributions on a non-isotropic 
pseudo-Riemannian symmetric space 
of rank one}%
\author{Hiroyuki Ochiai
\thanks{The author is supported in part by a Grant-in-Aid for 
Scientific Research (B) 15340005 from the Ministry of Education,
Culture, Sports, Science and Technology.}
}
\date{Dedicated to the sixty-fifth birthday\\
of Professor dr. Gerrit van Dijk}
\newcommand{\address}{\noindent
E-mail:  ochiai@math.nagoya-u.ac.jp \\
Department of Mathematics, Nagoya University, \\
Chikusa, Nagoya 464-8602, Japan. 
}
\newcommand{\Abstract}{{\noindent{\bf{Abstract}}.  %
We investigate the structure of invariant distributions %
on a non-isotropic non-Riemannian symmetric space of rank one. %
Especially, the $J$-criterion related to the generalized Gelfand pair %
is shown for this space without imposing the condition %
on the eigenfuction of the Laplace-Bertrami operator.  %
}}
\usepackage{amsmath}
\usepackage{amssymb}
\newtheorem{proposition}{{\bf Proposition}}
\newtheorem{theorem}[proposition]{{\bf Theorem}}

\newtheorem{lemma}[proposition]{{\bf Lemma}}
\newtheorem{remark}[proposition]{{\bf Remark}}

\newcommand{\qed}{\hspace*{\fill} \fbox{\;}\par\bigskip}

\newcommand{\Z}{{\mathbf Z}}
\newcommand{\R}{{\mathbf R}}

\newcommand{\diag}{{\mathrm{diag}}}
\newcommand{\trans}{{{}^t}}
\newcommand{\tG}{{\tilde G}}
\newcommand{\tH}{{\tilde H}}

\begin{document}
\maketitle

\Abstract

\section{Introduction}


In this paper,
the problem \cite{Thomas4} to understand the double coset space
$H \backslash G /H$
in terms of functions of distribution class
for a homogeneous space $G/H$
is examined in
one specific case study
for the space $G/H=SL(n+1,\R)/GL^+(n,\R)$.
If one deals with the class of arbitrary functions,
then 
the $H \times H$-orbit decomposition on $G$
gives enough information.
This is the case for the symmetric spaces over the finite field
(c.f. \cite{vD95})
and the set of characteristic functions of orbits forms a linear basis of 
invariant functions. 
If one deals with the regular functions on the algebraic group,
it is enough to look at the ring of invariants.
More generally,
if one considers bi $H$-invariant continuous functions,
then the lower-dimensional orbits can be ignored.

We here consider distributions.
The group case,
the famous result of Harish-Chandra tells us
that every invariant eigendistribution is locally integrable
and that there is no singular invariant eigendistribution.
For  Riemannian symmetric spaces, 
every invariant eigendistribution is real analytic.
However, for the case of non-Riemannian symmetric spaces,
the natural function space for invariant eigendistributions
(spherical functions)
is the class of distributions
(or its generalization to hyperfunctions).
A lot of works have been done for spherical functions
and also the eigenspace representations,
which imposes the eigenspace condition without invariance.
Here we will consider invariant distributions
without the eigenfunction condition of (indefinite) Laplacian.
In the case of invariant eigendistributions,
the space is finite-dimensional,
and especially for the rank-one case
we can write down explicitly such distributions  
using e.g., hypergeometric functions.
On the other hand, the space of invariant distributions
is infinite-dimensional, and has no ring structure in general.

Our target space $G/H$ is a symmetric space,
however,
we use the related homogeneous spaces,
having a larger symmetry.
In the case of the analysis on the tangent space of the symmetric space, 
such an enlargement of the symmetry, which has been discussed in
\cite{Kowata} (see also \cite{O}), can be formulated on the tangent space itself,
but on the symmetric spaces
it is necessary to introduce such homogeneous spaces
to formulate the enlargement of the symmetry.

Anyway, for the symmetric space or its related (spherical)
homogeneous space,
the double coset space $H \backslash G/H$
is not a manifold especially near the origin.
It is neither smooth nor Hausdorff, in general.
Since we are considering distributions,
the analysis at the origin is complicated on such a space.
For a non-archimedian ($p$-adic) local field,
the structure of distributions
are easier than for real numbers $\R$,
e.g., the extensions and decomposition of distributions reduces 
some problem to the orbit decompositions,
and the structure of the distributions supported on a lower-dimensional
submanifold is simpler (there is no derivative of Dirac delta).
c.f. \cite{BovD}.
In our Lie group case, analysis around some singular locus
will be more subtle, and it is a point of discussion.
We apply the result \cite{O} on the tangent space 
of $G/H$ to prove the extra symmetry which
every invariant distribution has.
This enables us to eliminate a contribution of some singular loci
on the non-isotropic symmetric spaces.
This is given in Section~\ref{sec:inv}.
In Section~\ref{sec:geo},
we summarize several geometric facts,
which will be used in the later sections
or help us to understand the double coset space.

As an application of the main theorem in Section~\ref{sec:inv},
we discuss the property called {\it generalized Gelfand pair}
in Section~\ref{sec}.
This notion is related to the uniqueness
of the decomposition of the left regular representations
of $L^2(G/H)$,
see \cite{vD93}.
So-called $J$-criterion is
one of the well-known sufficient conditions to be a generalized Gelfand pair.
The $J$-criterion for a class of semisimple symmetric spaces
of rank one is examined in \cite{vD86}.
Here we slightly generalize the statement
by dropping the eigenfunction condition.
Note that the $J$-criterion holds while the geometric counterpart
does not hold. 
It suggests a careful study between 
the orbit spaces and the space of invariant distributions.

The author would like to his gratitude to Professor dr. van Dijk,
and organizing committee of the conference.

\section{Geometric background}
\label{sec:geo}
\subsection{Enlargement of the symmetry of the space}

Let $G=SL(n+1,\R)$,
and $\sigma$ be the involution of $G$  defined by the conjugation by 
the matrix $\diag(-1, 1,\dots,1)$.
The fixed point subgroup is
\[
G^\sigma=\left\{\left(\begin{array}{cc}
(\det h)^{-1} & 0 \\ 0 & h \end{array}\right) \mid h \in GL(n,\R) \right\},
\]
which is often to be identified with $GL(n,\R)$.
The identity component of $G^\sigma$ is denoted by $H$,
which is isomorphic to $GL^+(n, \R)$.
The homogeneous space $G/H$ is
a non-isotropic non-Riemannian symmetric space of rank one,
which is our main concern in this paper.
We introduce $H_1 =SL(n,\R) \subset H=GL^+(n,\R)$,
then we have $G^\sigma/H=\Z/2\Z$,
$G^\sigma/H_1\cong GL(1,\R) =\R^\times$,
$H/H_1 \cong GL^+(1,\R)$, and $H \cong H_1 \times GL^+(1,\R)$.
Then we have $G/H=(G/H_1) / (H/H_1) = (G/H_1) / GL^+(1,\R)$.
It is easy to see that
the normalizer of $H_1$ in $G$ is $G^\sigma$.

We will describe the homogeneous space $G/H_1$ as follows.
Let
$X_1 = \{ (x,y) \mid x,y \in \R^{n+1}, \langle x, y \rangle = 1\}$.
Here we regard $x=\trans(x_1,x_2,\dots, x_{n+1}) = \trans(x_1,x')$,
$y=\trans(y_1,y_2,\dots,y_{n+1})=\trans(y_1,y')$ as column vectors.
The action of $G$ on $X_1$ is defined by
\[
g. (x,y) = (g x, \trans g^{-1} y),
\qquad x, y \in \R^{n+1}, g \in G.
\]
It is easy to see that this action is transitive.
We set $x_1=\trans (\trans e_1,\trans e_1)$ with $e_1=\trans(1,0,\dots,0)\in \R^{n+1}$.
Then the isotropy subgoup of $G$ at $x_1$ is $H_1$,
and we have a natural isomorphism $X_1 =G/H_1$.

We define the function $Q_0(x,y) = x_1y_1+x_2y_2+\cdots+x_{n+1}y_{n+1}$
on $\R^{2(n+1)}$.
We denote by $\tG=SO_0(Q_0)\cong SO_0(n+1,n+1)$ the identity component
of the orthogonal group $O(Q_0)$ corresponding to the quadratic form $Q_0$.
Then $\tG$ also acts on $X_1$ transitively.
The isotropy subgroup $H_2$ of $\tG$ at $x_1$ is isomorphic to
$SO_0(n,n+1)$.
Then $X_1 \cong G/H_1 \cong \tG/H_2$.
The role of this isomorphism for the harmonic analysis has been emphasized 
in \cite{Ko}, Example 5.2.
This isomorphism means that the space $G/H_1$ has a larger symmetry $\tG$.
The expression $X_1=G/H_1$ is not a symmetric space but
a (real form of) spherical homogeneous space.
The expression $X_1=\tG/H_2$ is an isotropic symmetric space of rank one.

\subsection{Invariants and orbits}
The function
$Q(x,y)=1-x_1 y_1 =  \trans x' y'$ on $X_1$
is $H$-invariant.
The functions $x_1$ and $y_1$ on $X_1$ are 
$H_1$-invariant.
The map
\[
q: X_1 \ni (x,y) \mapsto (x_1,y_1) \in \R^2
\]
is real-analytic, surjective, and 
$H_1$-invariant.
As is seen later, 
the map $q$ almost classifies the set of $H_1$-orbits on $X_1$,
and any $H_1$-invariant continuous function on $X_1$ 
is a pull back by the map $q$ of a
continuous function on $\R^2$.
For distributions, the question is more subtle.

Let
$\tH=SO_0(Q_1)$ be the identity component of the orthogonal group
$O(Q_1)$ corresponding to the quadratic form 
$Q_1(x',y')=x_2y_2+\cdots+x_{n+1}y_{n+1}$ on $\R^{2n}$.
Then $\tH$ is identified with a subgroup of $\tG=SO_0(Q_0)$.
The group $\tH$ preserves the invariants $x_1$ and $y_1$.
Now we give several orbit decompositions.

(1) $H_1$-orbit decomposition on $X_1$.
The map $q=(x_1,y_1)$ almost classifies the $H_1$-orbits on $X_1$.
Say, for $(x_1,y_1) \in \R^2$ with $x_1 y_1 \neq 1$,
the fiber $q^{-1}(x_1,y_1)$ is an $H_1$-orbit
if  $n\ge2$.
For $x_1 \neq 0$,
the fiber $q^{-1}(x_1,x_1^{-1})$
splits into four $H_1$-orbits
$\{ (x',y') \mid \trans x' y' = 0, x' \neq 0, y'\neq 0\}$,
$\{ (0, y') \mid y' \neq 0 \}$,
$\{ (x', 0) \mid x' \neq 0 \}$,
and $\{ (0,0) \}$
if $n\ge3$.
If $n=2$, then
the fiber $q^{-1}(x_1,x_1^{-1})$
consists of $H_1$-orbits
$\{ (x',y')=(x_2,x_3,y_2,y_3) \mid x' \neq 0, y_2=-t x_3, y_3 = t x_2 \}$
with $t \neq 0$,
$\{ (0, y') \mid y' \neq 0 \}$,
$\{ (x', 0) \mid x' \neq 0 \}$,
and $\{ (0,0) \}$.

(2) $GL^+(n,\R)$-orbit decomposition on $X_1$.
We consider the action of $GL^+(n,\R)$ on  $X_1$ by
\[
(x_1,x',y_1,y') \mapsto (x_1, hx', y_1, \trans h^{-1} y')
\mbox{ for }
h \in GL^+(n,\R), (x,y) \in X_1.
\]
The orbit decomposition on $X_1$ under $GL^+(n,\R)$
is the same as that under $H_1$, except for
the case $n=2$ and the fiber $q^{-1}(x_1,x_1^{-1})$,
which is decomposed into five $GL^+(n,\R)$-orbits
$\{ (x',y')=(x_2,x_3,y_2,y_3) \mid x' \neq 0, \pm  (x_2 y_3-x_3 y_2)>0 \}$,
$\{ (0, y') \mid y' \neq 0 \}$,
$\{ (x', 0) \mid x' \neq 0 \}$,
and $\{ (0,0) \}$.

(3) $\tH$-orbit decomposition on $X_1$.
Since the map is $\tH$-invariant, 
the fiber $q^{-1}(x_1,y_1)$
for $(x_1,y_1) \in \R^2$ with $x_1 y_1 \neq 1$ is an $\tH$-orbit
if  $n\ge2$.
The fiber $q^{-1}(x_1,x_1^{-1})$
consists of two $\tH$-orbits,
$\{ 0 \}$ and $\{(x',y') \neq 0 \mid \trans x' y' = 0 \}$.

(4) $H$-orbit decomposition on $X$.
This is equivalent to $GL^+(n,\R) \times GL^+(1,\R)$-orbits
on $X_1$.
In this case, the map $Q : X_1 \ni (x,y) \mapsto 1-x_1y_1 \in \R$
almost classifies orbits.
In fact, for $t \neq 0,1$,
the fiber $Q^{-1}(t)$ 
consists of two orbits
$\{(x,y) \in X_1 \mid x_1 y_1 = 1-t, \pm x_1>0 \}$
if  $n\ge2$.
$Q^{-1}(1)$ consists of five orbits
$\{(x,y) \in X_1 \mid \pm x_1>0 , y_1=0 \}$,
$\{(x,y) \in X_1 \mid x_1 =0, \pm y_1>0 \}$,
and $\{(x,y) \in X_1 \mid x_1=y_1=0\}$
if $n \ge 2$.
For $x_1 \neq 0$,
the fiber $Q^{-1}(0)$
splits into eight orbits
$\{ (x_1,x',x_1^{-1}, y') \mid \pm x_1>0, \trans x' y' = 0, x' \neq 0, y'\neq 0\}$,
$\{ (x_1,0,x_1^{-1}, y') \mid \pm x_1>0, y' \neq 0 \}$,
$\{ (x_1,x',x_1^{-1}, 0) \mid \pm x_1>0,x' \neq 0 \}$,
and $\{ (x_1,0,x_1^{-1},0) \mid \pm x_1>0 \}$
if $n\ge3$.
If $n=2$, then
the fiber $q^{-1}(0)$
consists of ten orbits
$\{ (x_1,x_2,x_3,x_1^{-1}, y_2,y_3) \mid 
\pm x_1>0, x' \neq 0, \pm  (x_2 y_3-x_3 y_2)>0 \}$,
$\{ (x_1,0, x_1^{-1}, y') \mid \pm x_1>0, y' \neq 0 \}$,
$\{ (x_1,x', x_1^{-1}, 0) \mid \pm x_1>0, x' \neq 0 \}$,
and $\{ (x_1,0,x_1^{-1},0) \mid \pm x_1>0 \}$.

(5) $\tH \times GL^+(1,\R)$-orbit decomposition on $X_1$.
For $n\ge 2$ and $t \neq 0$,
the orbit decomposition of the fiber $Q^{-1}(t)$
is the same as that of (4).
The fiber $Q^{-1}(0)$
splits into four orbits
$\{ (x_1,x',x_1^{-1}, y') \mid \pm x_1>0, \trans x' y' = 0, (x',y') \neq 0 \}$,
and $\{ (x_1,0,x_1^{-1},0) \mid \pm x_1>0 \}$
if $n\ge3$.
If $n=2$, then
the fiber $Q^{-1}(0)$
consists of four orbits
$\{ (x_1,x_2,x_3,x_1^{-1}, y_2,y_3) \mid 
\pm x_1>0, (x',y') \neq 0 \}$,
and $\{ (x_1,0,x_1^{-1},0) \mid \pm x_1>0 \}$.

The geometric $J$-criterion (for an involution $\theta$) is
the statement that $H g^{-1} H = H \theta(g) H$ for all $g \in G$.
Note that for the involution $\theta(g)=\trans g^{-1}$,
the space $G/H$ does not satisfy the geometric $J$-criterion
since the $H$-orbit 
$\{ (x_1,0,x_1^{-1}, y') \mid \pm x_1>0, y' \neq 0 \}$
is mapped to
$\{ (x_1,x',x_1^{-1}, 0) \mid \pm x_1>0,x' \neq 0 \}$,
respectively,
by $HgH \mapsto H\theta(g)^{-1}H$.
Nevertheless,
we will prove the (original) $J$-criterion (for distributions)
in Section~\ref{sec}.

\section{Bi-invariant distributions}
\label{sec:inv}

We denote by $C^{-\infty}$ the set of functions of distribution class.
We simply call a function of distribution class a {\it distribution}.

We start from a direct conclusion of the enlargement of the symmetry
of the space $G/H_1$.
\begin{lemma}
There  are natural identifications
between the set of distributions with the following properties.
\begin{enumerate}
\item[(i)]
An $H$-bi-invariant distribution on $G$.
\item[(ii)]
A left $H$ and right $H/H_1$-invariant distribution on $G/H_1$.  
\item[(iii)]
A left $GL^+(n,\R) \times GL^+(1,\R)$-invariant distribution on $X_1$.
\end{enumerate}
\end{lemma}
Proof. 
Since  $H$ normalizes $H_1$,
$H$ acts on  $G/H_1$ from the right as
\[
G/H_1 \ni g H_1 \mapsto g h H_1 \in G/H_1,
\]
which induces the action of $H/H_1\cong GL^+(1,\R)$ on $X_1$
as $(x,y) \mapsto ((\det h)^{-1} x, (\det h) y)$,
that is, $(x,y) \mapsto (t^{-1} x, t y)$
by $t>0$.
On the other hand,
the action of $H$ on $G/H_1$ from the left is
\[
(x_1,x',y_1,y') \mapsto
((\det h)^{-1} x_1, hx', (\det h)y_1, \trans h^{-1} y').
\]
Then, an $H$-bi-invariant function on $G$
is identified with  a function on $X_1$ invariant under
the action
\[
(x,y) \mapsto (t x, t^{-1}y),
(x_1,x',y_1,y') \mapsto (x_1, h x', y_1, \trans h^{-1}y'),
\]
for all $t>0$ and $h \in GL^+(n,\R)$.
Then $GL^+(1,\R) \times GL^+(n,\R)$
is a subgroup of $\tG$.
\qed
It is somewhat mysterious that
the left-right action of $H$ on $G$ 
turns to be equivalent to the left action of 
$H \times GL^+(1,\R)\subset \tG$.

The next theorem
shows the enlargement of the symmetry on invariant distributions,
which has not been predicted by the geometry,
e.g., orbit structures.

\begin{theorem}
Let $n \ge 3$.
Then 
a $GL^+(n,\R)$-invariant distribution on $X_1$ is $\tH$-invariant.
\end{theorem}
Proof : 
The invariance is local, so we may consider the restrictions
on the open subset $\{Q<1\}$ and $\{Q>0\}$.
(i)
On the open subset $Q<1$.
We can take the coordinates $(x_1,x',y')$ on 
the open subset $\{(x,y) \in X_1 \mid Q<1, x_1 \neq 0\}
\cong \{(x',y') \in \R^{2n} \mid \trans x' y' <1 \} \times \{ x_1 \neq 0 \}$.
Actually, $y_1 = (1-\trans x' y')/x_1$.
The action of $\tH$ is only on the variables $(x',y')$.
It has been proved in  Theorem~1 of \cite{O}  that
the D-modules for $H$-invariants equals to that for $\tH_1$.
It means that any $H$-invariant distribution on this space
is $\tH$-invariant.
On the open set $y_1 \neq 0$, we have the same argument.

(ii)
On the open subset $Q>0$.
We can take the coordinates 
$\{(x,y) \in X_1 \mid Q>0 \}
\ni (x_1,x',y_1,y') \mapsto
((x', y'/(1-x_1 y_1)), (x_1,y_1)) \in
\{ (\xi, \eta) \in \R^{2n} \mid \trans \xi \eta =1 \} 
\times \{ (x_1,y_1) \mid x_1 y_1 <1 \}$.
Then both the groups $H$ and $\tH$ act on the first factor
transitively.  So, the $H$-invariants implies $\tH$-invariants.
\qed
\begin{remark}
(1)
The proof shows that the same statement for the theorem holds
if we replace `distribution' by `hyperfunction'.
(2)
We can regards Lemma~7.4 of \cite{KovD}
as a special case of theorem;
they have proved that an $H$-invariant distribution on $X$
supported on $\{(x,y)\in X_1 \mid x'=0 \mbox{ or } y'=0\}$
is supported on $\{(x,y) \in X_1 \mid x'=0 \mbox{ and } y'=0\}$.
The proof uses the fact that distributions are of finite order.
\end{remark}

In the case of the tangent spaces,
such an extension of the symmetry has
been observed in \cite{Kowata} with the eigenfunction
condition, and is extended in \cite{O}
without the eigenfunction condition. 
These works are inspired by  \cite{vD84}.
Compared to the case of the tangent spaces,
the geometric setting is slightly more subtle
for the case of homogeneous spaces.

\section{Generalized Gelfand pairs}
\label{sec}

We define an anti-involution $J:G \rightarrow G$ by $J(g)=\trans g$ for $g \in G$.
The map $J$ induces the linear endomorphism $J^*$ on
$C^{-\infty}(G)^{H \times H}$ by $f \mapsto f \circ J$.
The natural isomorphism $i : G/H_1 \rightarrow X_1$
induces the isomorphism $i^* : C^{-\infty}(X_1)^{H \times GL^+(1,\R)}
\rightarrow C^{-\infty}(G)^{H \times H}$.
We now describe the map $J^*$ on
$C^{-\infty}(G)^{H \times H} \cong C^{-\infty}(X_1)^{H \times GL^+(1,\R)}$.
We define $j : X_1 \ni (x_1,x',y_1,y') \mapsto (x_1,-y',y_1,-x') \in X_1$.
\begin{lemma}
For any $f \in C^{-\infty}(X_1)^{H \times GL^+(1,\R)}$,
we have
$J^* i^* f = i^* j^* f$.
\end{lemma}
Proof:
(i)
On $Q<1$.
Let $U=\{(x,y) \in X_1 \mid x_1y_1>0\}$,
\[
\alpha : U \times H_1 \ni ((x,y), h) \mapsto
\left(\begin{array}{ccc}
x_1 & - y_2 & - \trans y''/y_1 \\
x_2 & y_1 & 0 \\
x'' & 0 & I_{n-1}
\end{array}
\right) \in G,
\]
where $x'=\trans(x_2,x'')$, $y'= \trans(y_2, \trans y'')$
and
$\varphi : U \times H_1 \rightarrow GL^+(n,\R) \subset \tH$
by
$\varphi((x,y),h)=\diag(y_1,1,\dots,1)/y_1$.
The image of $\alpha$ is an open subset $i^{-1}(\{Q>0\})$ of $G$.
We have
$i \circ J \circ \alpha 
= \varphi \cdot (j \circ i \circ \alpha)$.
This proves
$
(J^* i^* f)(\alpha((x,y),h))
= f(i \circ J \circ \alpha((x,y),h)))
= f(\varphi((x,y),h) j\circ i \circ \alpha((x,y),h))
= f(j\circ i \circ \alpha((x,y),h))
= (i^* j^* f)(\alpha((x,y),h))$.

(ii) 
On $Q>0$.
Let $U=\{(x_1,y_1) \in \R^2 \mid x_1y_1<1 \}$,
\[
\alpha : U \times H  \times H \ni ((x_1,y_1),h,h')
\mapsto h \left(\begin{array}{ccc}
x_1 & -1+x_1y_1 & 0 \\
1 & y_1 & 0 \\
0 & 0 & I_{n-1}
\end{array}\right) h' \in G,
\]
and 
$\varphi: U \times H \times H \rightarrow GL^+(n,\R) \subset \tH$
by
$\varphi(h,h')=\det(h h')^{-1} \trans(h h')
\in GL^+(n, \R)$.
The image of $\alpha$ is $i^{-1}(\{Q<1\})$.
We have
$i \circ J \circ \alpha 
= \varphi \cdot (j \circ i \circ \alpha)$.
This proves
$
(J^* i^* f)(\alpha((x_1,y_1),h,h'))
=  (i^* j^* f)(\alpha((x_1,y_1),h,h'))$.
\qed

\begin{theorem}
For $n \ge 3$,
$J^*$ is the identity on $C^{-\infty}(G)^{H \times H}$.
\end{theorem}
Proof.
It is enough to prove that $j^*$ is the identity
on $C^{-\infty}(X_1)^{\tH \times GL^+(1,\R)}$.
For any $f \in C^{-\infty}(X_1)^{\tH \times GL^+(1,\R)}$,
the support of $j^*f-f$ is contained in $\{(x,y) \in X_1 \mid x'=y'=0\}$.
Then the distributuion $j^* f-f$ can be written
as $p(\square) \delta(x',y')$ with some polynomial $p$
of the indefinite Laplacian 
$\square=\sum_{i=2}^{n+1} \partial^2/\partial x_i \partial y_i$
on each open subset $\{(x,y) \in X_1 \mid x_1 y_1>0, \pm x_1>0\}$.
This means that $(j^* f - f)$ is invariant under the action of $j$,
and that it is zero.
\qed

For the space $G/H$,
the $J$-criterion 
that any bi-$H$-invariant {\it eigen}distribution on
$G$ is invariant under $J^*$
has been proved in \cite{vD86}.
The $J$-criterion implies 
that the space $G/H$ is a generalized Gelfand pair
\cite{Thomas}, \cite{vD86}.


\address

\end{document}